\input amstex
\input epsf
\documentstyle{amsppt}
\NoBlackBoxes 
\hfuzz 6pt

\hoffset=1.5cm
\hsize=5.0truein
\vsize=8.05truein

\def\C{\Bbb C}
\def\R{\Bbb R}
\def\H{\Cal H}
\def\fy{\varphi}
\def\ls{\leqslant}
\def\gs{\geqslant}

\baselineskip=10pt
\line{\eightrm Journal of Knot Theory and Its Ramifications\hfil}
\line{\eightpoint $\copyright$\, World Scientific Publishing
Company\hfil}
\vglue 5pc
\baselineskip=13pt

\nologo

\topmatter
\title
Topological classification of generic real rational functions 
\endtitle
\author Sergei Natanzon$^{*,}$\footnote 
{\eightrm
Partly supported by the INTAS grant 00-0259 and by the RFFI grant
01-01-00739.\kern 1.5cm},
Boris Shapiro$^{\dag}$ and Alek Vainshtein$^\ddag$ 
\endauthor
\affil 
$^*$ \eightit Moscow State University and Independent University of Moscow, Russia\\
$^{\dag}$ Department of Mathematics, University of Stockholm, Stockholm, 
Sweden\\
$^\ddag$ Department of Mathematics and Department of Computer Science\\ 
University of Haifa, Haifa, Israel 31905.
\endaffil
\endtopmatter

\nopagenumbers
\headline={\ifodd\pageno\rightheadline \else \leftheadline\fi}
\def\rightheadline{\hfil\eightit  Classification of Real Rational Functions
\quad\eightrm\qquad}
\def\leftheadline{\eightrm\qquad\quad \eightit  S.~Natanzon, 
B.~Shapiro \& A.~Vainshtein\hfil}

\document
\vglue 10pt
\vglue 16pt
\centerline{\eightrm ABSTRACT}

{\rightskip=1.5pc
 \leftskip=1.5pc
 \eightrm\baselineskip=10pt\parindent=1pc
To any real rational function with generic ramification points we assign
a combinatorial object, called a garden, which consists of a weighted 
labeled directed planar chord diagram and of a set of weighted rooted trees
each corresponding to a face of the diagram. We prove that any garden 
corresponds to a generic real rational function, and that equivalent functions
have equivalent gardens. 
\vglue 10pt	
\noindent
{\eightit Keywords}\/: Real rational functions; Generic ramification;
Topological invariants; Hurwitz numbers.
\vglue 18pt}

\baselineskip=13pt	

\line{\bf  1. Introduction \hfil}
\vglue 5pt

Let $f\: P\to\bar\C=\C\cup\infty$ be a meromorphic function of degree $n$
on a compact Riemann surface $P$ of genus $g$. We say that $f$ is a
{\it generic\/} (complex) meromorphic function if the preimage $f^{-1}(z)$
of any point $z\in\bar\C$ consists of either $n$ or $n-1$ points; 
equivalently, the singularities of $f$
are of degree two, and at any two distinct singular points $f$  takes distinct
values. The points $z$ for which  $|f^{-1}(z)|=n-1$
are called {\it simple ramification points}. The set of
all simple ramification points of $f$ is denoted $\Sigma(f)$; 
by the Riemann--Hurwitz formula, it consists 
of $2n+2g-2$ points.

Two meromorphic functions $f_i\: P_i\to\bar\C$ ($i=1,2$) are called {\it 
equivalent\/} if there exists a biholomorphic map $\fy\: P_1\to P_2$ such
that $f_1=f_2\circ \fy$. Let $\C \H_{g,n}$ be the set of equivalence classes of
complex generic meromorphic functions of degree $n$ on surfaces of genus 
$g$.
The correspondence $f\mapsto \Sigma(f)$ generates a covering 
$\C\Phi_{g,n}\:\C \H_{g,n}\to \C Q_{g,n}$, where $\C Q_{g,n}$ is the
configuration space of all $(2n+2g-2)$-tuples of unordered distinct points on 
$\bar\C$,
or, equivalently, the projectivized space of complex homogeneous degree
$2n+2g-2$ polynomials in two variables without multiple roots.
We assume that $\C \H_{g,n}$ is provided with the weakest topology for which
the map $\C\Phi_{g,n}$ is continuous. According to the
Hurwitz theorem \cite{Hu}, $\C \H_{g,n}$ is a connected space. 
 The degree 
$\C h_{g,n}$ of the covering $\C\Phi_{g,n}$, and its analogs for arbitrary
meromorphic functions, are called the {\it Hurwitz numbers}. These numbers
arise in many situations in mathematical physics; for example, they
generate correlators of topological field theory, see \cite{CMR}. 
In recent years
they attracted much attention. For the case $g=0$, the Hurwitz
numbers are, in particular, calculated in \cite{CT}: 
$$
\C h_{0,n}=\frac{n^{n-3}(2n-2)!}{n!};
$$
in fact, this result was apparently known already to Hurwitz himself. For the 
case $g=1$, the Hurwitz numbers are calculated in \cite{GJV}:
$$
\C h_{1,n}=\frac1{24}\left(n^n-n^{n-1}-\sum_{i=2}^n\binom{n}{i}(i-2)!
n^{n-i}\right).
$$
For certain classes of nongeneric rational functions, Hurwitz numbers
where studied in \cite{SSV, GL, ELSV}. The former paper exploits the classic
approach due to Hurwitz, which links the numbers in question to the characters
of the symmetric group. The approach developed in the other two papers
is due to Arnold and is based on the singularity theory.

In the present paper we study a similar problem for real meromorphic
functions. A real meromorphic function is defined on a {\it real algebraic 
curve}, which is a pair $(P,\tau)$, where $P$ is a complex algebraic curve
(a compact Riemann surface), and $\tau\:P\to P$ is the antiholomorphic 
involution (the involution of complex conjugation). A {\it real meromorphic 
function\/} is a complex meromorphic function $f\: P\to\bar \C$ such that
$\overline{f(\tau p)}=f(p)$ for any $p\in P$, see e.g. \cite{N2}. A real 
meromorphic function $(P,\tau,f)$ is said to be {\it generic\/} if $(P,f)$
is a generic complex meromorphic function. Evidently,
for any real meromorphic function $f$ one has $\overline{\Sigma(f)}=\Sigma(f)$.
Two real meromorphic functions
$(P_i,\tau_i,f_i)$ ($i=1,2$) are called {\it equivalent\/} if there exists
a biholomorphic map $\fy\: P_1\to P_2$ such that $f_1=f_2\circ \fy$ and
$\fy\circ \tau_1=\tau_2\circ \fy$. Let $\R \H_{g,n}$ denote the space of 
equivalence
classes of generic real meromorphic functions of degree $n$ on surfaces 
of genus $g$.
The topology of $\C \H_{g,n}$ generates a topology on $\R \H_{g,n}$; in this
topology $\R \H_{g,n}$ is not connected (see \cite{N2}). 
The covering
$\C\Phi_{g,n}$ generates a covering $\R\Phi_{g,n}\:\R \H_{g,n}\to\R Q_{g,n}$,
where $\R Q_{g,n}$ is the projectivized space of real homogeneous degree 
$2n+2g-2$
polynomials in two variables without multiple roots, see \cite{N2}.

In this paper we study connected components of $\R \H_{0,n}$; the points
of this space are called equivalence classes of
{\it generic real rational functions}. Allowing a slight abuse of language,
we refer to the elements of $\R\H_{0,n}$ as generic real rational
functions, and write $f\in\R\H_{0,n}$ meaning that the equivalence class
of $f$ belongs to $\R\H_{0,n}$. We define
topological invariants that distinguish each connected component $H\subset
\R \H_{0,n}$ and find the corresponding Hurwitz number $\R h_H$, that is, 
the degree of the restriction of $\R\Phi_{0,n}$ to $H$. 
For certain classes of nongeneric real rational functions Hurwitz numbers
were studied in \cite{Ar, Ba, Sh, SV1}. Note that interesting cell 
decompositions of the space of complex generic rational 
(and, more generally, meromorphic) functions were studied in \cite{CP, BC}. 

The principal result of this note is as follows; see precise definitions in 
\S 2. 

\proclaim{Main Theorem} The set of all connected components of the space
 $\R \H_{0,n}$ is in a 1-1-correspondence with the set of the equivalence
classes of all gardens of weight $n$.
\endproclaim

The structure of the paper is as follows. In \S 2 we introduce the notion 
of a garden and construct it for any given generic real rational function. 
In \S 3 we prove the main theorem as well as all necessary preliminary 
results 
for the calculation of Hurwitz numbers, which is carried out in \S 4. 
Finally, \S 5 contains some open questions and comments. 

The authors are sincerely grateful to the Max-Planck Mathematical Institute 
 in Bonn for the financial support and excellent research atmosphere during 
 the fall of 2000 when this project was started. The second author wants 
 to acknowledge the hospitality of IHES, Paris in January 2001 during the 
 final stage of preparation of the manuscript.

\vglue 12pt
\line{\bf 2. Topological Invariants of Generic Real Rational Functions\hfil}
\vglue 5pt

The purpose of this section is to introduce a combinatorial object, which we 
assign to any
generic real rational function. This object, called a {\it garden}, consists
of a weighted labeled directed planar chord diagram and of a set of
weighted rooted trees each corresponding to a face of the diagram.

\subheading{2.1. Defining gardens abstractly}
By a {\it planar chord diagram\/} (of {\it order\/} $2l$) we mean a circle 
drawn on the plane together with
$2l$ points on this circle partitioned into $l$ pairs in such a way that for
any two pairs, the chords joining the points from the same pair do
not intersect. The above $2l$ points are called the {\it vertices\/} of the 
chord diagram; the chords joining the vertices from the same pair, as well as 
the arcs of the circle joining adjacent vertices, are called the {\it edges}.
Clearly, a planar chord diagram is a plane graph, so the notion of its
{\it faces\/} is defined in a usual way (except for the outer face of the 
graph, which is not a face of the diagram). We say that a planar chord
diagram is {\it directed\/} if its edges are directed in such a way that
the boundary of each face becomes a directed cycle. Obviously, in order 
to direct a planar chord diagram it suffices to direct any one of its edges. 
Therefore there exist exactly two possible ways of directing a diagram, 
which are 
opposite to each other, i.e. the second one is obtained from the first one by
reversing the direction of every edge. Once and for all 
fixing the standard orientation of the plane, we call a face of a 
directed planar chord diagram {\it positive\/} if the face lies to the left 
when we traverse its boundary according to the chosen direction, and {\it 
negative} otherwise. All neighbors of positive faces are negative, and vice
 verse. 
 
A planar chord diagram is said to be {\it weighted\/} if each edge is 
equipped with a nonnegative integer weight, and {\it labeled\/} if there exists
a bijection $\beta$ ({\it labeling\/}) that takes the vertex set of the
diagram to the set $\{1,2,\dots,2l\}$. Two labelings $\beta_1$ and $\beta_2$ 
are said to be {\it cyclically equivalent\/} if
$\beta_1(v)-\beta_2(v) \mod 2l$ is a constant not depending on the choice of a
vertex $v$.

Consider a labeled directed planar chord diagram. For any face $j$ we denote
by $d_j$ the number of descents in the sequence of vertex labels ordered
cyclically along the boundary of the face; clearly $d_j\gs1$. 
If the diagram is also weighted,
we denote by $t_j$ the sum of $d_j$ and the weights of all the edges along
the boundary of the face $j$. 

Recall that a {\it rooted tree\/} is a tree with one distinguished node 
called the {\it root\/}; all the other nodes of the tree are said to be
{\it inner}. Given a rooted tree with the root $r$ and two nodes
$u$ and $v$, we say that $u$ is a {\it child\/} of $v$ if the tree contains
the edge $(u,v)$, and if $v$ lies on the unique path between $u$ and $r$.
We say that a rooted tree is {\it weighted\/} if each of its nodes is 
equipped with a positive integer weight.

\vskip 10pt
\centerline{\hbox{\epsfxsize=11cm\epsfbox{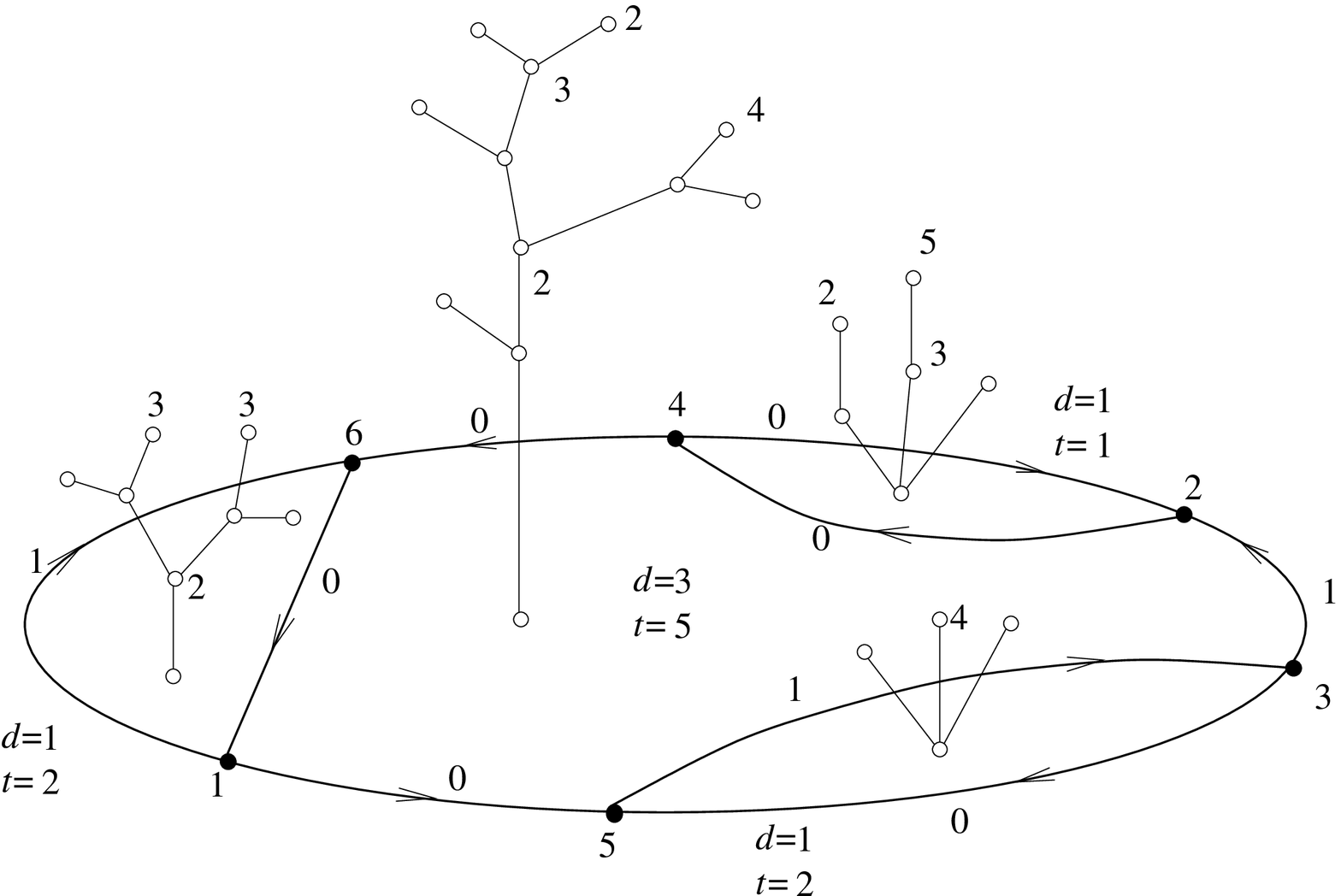}}}
\midspace{1mm}
\caption{Fig.~1. A garden of order 6 and total weight 106}
\vskip 5pt

{\smc Important definition.} A {\it garden\/} is a weighted labeled  
directed planar chord diagram
with a weighted rooted tree (possibly  consisting just of its root) 
corresponding to each face of the diagram. 
The weights of the {\it inner\/} nodes of 
the trees are arbitrary positive integers,
and the weight of the root of the tree corresponding to the face $j$ equals
$t_j$. The {\it total weight\/} of the garden equals twice the sum of the
weights of all the inner nodes of all trees plus the sum of the weights
of all roots.

An example of a garden is given on Fig.~1. The order of the diagram
equals $6$. The numbers written near the vertices and the edges are the
labels and the weights, respectively. The weights of the nodes are equal 
to one, unless specified otherwise. The total weight of the garden equals
$106$.

Two gardens are said to be {\it equivalent\/} if there exists a bijection
of the vertex sets of the corresponding chord diagrams that preserves
chords, their orientation, labels (up to the cyclic equivalence), 
rooted trees, and weights.

\subheading{2.2. Getting gardens from rational functions}  
To each function $f\in \R \H_{0,n}$ we associate the garden $G(f)$ as follows.
First of all, represent $\Sigma=\Sigma(f)$ as $\Sigma=\Sigma_R\cup\Sigma_I$,
where $\Sigma_R$ is the set of real critical values of $f$
(not necessary finite), and 
$\Sigma_I$ is the set of its non-real critical values.
Consider the preimage $S(f)$ of the real line $\bar\R =\R\cup\infty$ under 
$f$. Evidently, $S(f)$
contains $\bar\R$ and is invariant under the standard involution.
All the critical points of $f$ that correspond to 
critical values in $\Sigma_R$ are real as well. Indeed, if $x$
is a critical point with a real critical value, then $\bar x$ is a critical
point with the same critical value; therefore,
$x=\bar x$, since $f$ is generic. A similar argument
shows that the number of such critical points is even; we denote it 
$2l(\Sigma)$.

For each
critical point as above, $S(f)$ contains exactly four arcs incident to it.
Two of these arcs are the arcs of $\bar\R \subset S(f)$, while the other two
interchange under the standard involution; in particular, the other endpoints
of these two arcs coincide. Moreover, these arcs do not intersect outside
$\bar\R \subset S(f)$, since such an intersection point would be a 
critical point with a real critical value. Therefore, these arcs together
with $\bar\R\subset S(f)$ define a 2-dimensional cell complex on $\bar\C$.
The 2-cells of this complex are called the {\it faces\/} of $S(f)$.
Besides, $S(f)$ contains a number of closed curves called {\it ovals}.
For the same reasons as above, no two ovals intersect, and each oval lies 
entirely inside one face. Observe that each face lies entirely in one of the
two hemispheres $\bar\C\setminus\bar\R$; moreover, the image of a face under
the standard involution is a face as well, and all the ovals lying inside
the former face are mapped bijectively to the ovals lying inside the latter
face.

To construct $G(f)$ we start from a planar chord diagram of order $2l(\Sigma)$.
The vertices of the diagram correspond to the critical points with real 
critical values, and the chords correspond to the arcs of $S(f)$ lying in
the upper hemisphere; thus, the faces of the diagram correspond to the faces
of $S(f)$ lying in the upper hemisphere.
The orientation of the edges is induced by the orientation of $\bar\R$
in the image. To define the labeling
of the chord diagram, consider the natural linear order $<$ on $\Sigma_R$
(if $\infty$ belongs to $\Sigma_R$, we assume that it is the biggest critical
value).
The label of a critical point equals the number of the corresponding critical
value under this order. To define the weights, consider an arbitrary point 
$x\in \bar\R\setminus\Sigma_R$ and for any given arc (or oval) define $w(x)$ 
as the number of preimages of $x$ lying on this arc (oval). The weight of the 
arc (oval) is then
defined as the minimum of $w(x)$ over all $x\in \bar\R\setminus \Sigma_R$. 

To  construct the rooted tree corresponding to a given face we proceed 
inductively.
The root of the tree corresponds to the boundary of the face; the inner 
vertices correspond to the ovals contained in the face under consideration. 
If there are no inner ovals, the tree consists only of its root. Otherwise, 
given an oval, the subtree rooted at the corresponding inner vertex
contains exactly the vertices whose ovals lie inside the given oval.  
The weight of an inner vertex is equal to the weight of the corresponding
oval.
An example of a face and the corresponding rooted tree is given on Fig.~2.

Observe that if $g$ belongs to the same equivalence class of generic real
rational functions as $f$, then the garden constructed for $g$ coincides 
with the one constructed for $f$.
 
\vskip 10pt
\centerline{\hbox{\epsfxsize=10cm\epsfbox{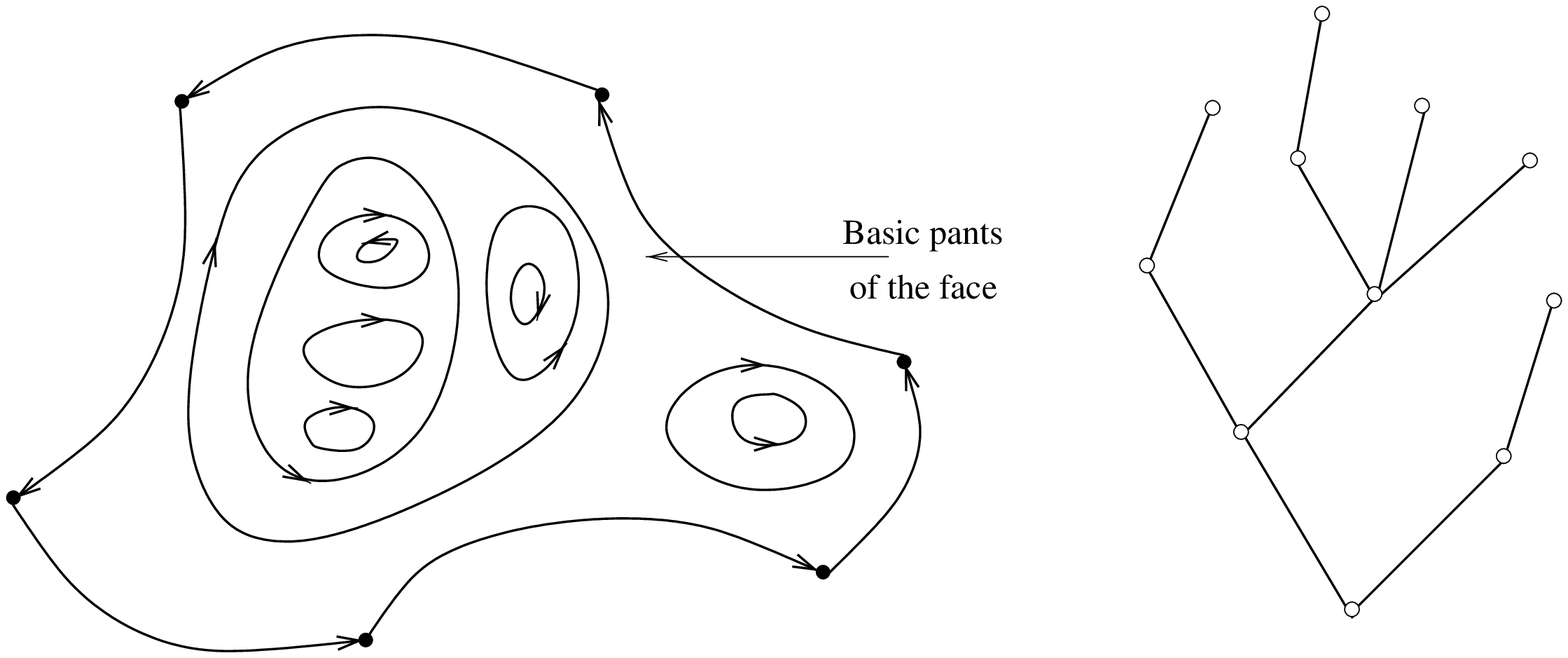}}}
\midspace{1mm}
\caption{Fig.~2. A face and the corresponding rooted tree}
\vskip 5pt

It is easy to see that the weight of a node coincides 
with the multiplicity of $f$ restricted to the corresponding oval
(or to the boundary of the corresponding  face). 
Since the total preimage of $\bar \R$ under a chosen $f\in \R \H_{0,n}$ 
coincides with $S(f)$, the total weight of its 
 garden $G(f)$ coincides with its degree and is therefore equal to $n$. 

 Given an abstract garden $G$ we can substitute each of its trees by the 
 appropriate system of weighted ovals, see  Fig.~2. Such a garden 
 will be called {\it represented\/}. In what follows we will freely 
 use both abstract and represented gardens. The connected components of the 
 complement to a represented garden are called {\it pants\/} 
 (as before, we disregard the outer face). Each pants 
 is a Riemann surface with a boundary consisting of a single outer 
 boundary component and some number of inner boundary 
 components. Note that we have assigned a certain weight to each connected 
 component of the boundary of each pants. Pants with weights on each 
 boundary component are called {\it weighted pants \/.} 
 The chosen direction of the chord diagram of a garden $G$ 
 extends in a unique way 
 to directions of all ovals such that every pants becomes either positive or 
 negative, i.e. lie either to the left (if the pants are positive) 
 or to the right (if the pants are negative) when we traverse 
 any component of the boundary of these pants.  The set of all weighted pants 
 of a given garden $G$ is called the {\it weighted pants collection\/} 
 and denoted by $\Pi(G)$.  

\vglue 12pt
\line{\bf 3. Realization Theorem and Connected Components of $\R \H_{0,n}$ 
\hfil}
\vglue 5pt

The Main Theorem is obviously equivalent to the following pair of statements. 

\proclaim{Theorem 1}  Let $\Sigma$ be an arbitrary set of $2n-2$ distinct
complex numbers invariant under the standard involution, of which exactly 
$2l$ are real. Any garden of order $2l$ and total weight $n$ is 
isomorphic to the garden $G(f)$ for some real meromorphic function  
$f\in\R \H_{0,n}$ such that $\Sigma(f)=\Sigma$.
\endproclaim

\medskip

\proclaim{Theorem 2}  Two rational functions belong to the same connected 
component of $\R \H_{0,n}$ if and only if they have equivalent gardens. 
\endproclaim

Both proofs require a number of additional statements. The idea of the proof
 of Theorem 1 is to construct a real topological covering 
 $\bar \C\to \bar \C$ with a given garden 
 and then, as usual in this field, to transform it into a holomorphic covering
  inducing the holomorphic structure on the preimage $\bar \C$ from that on 
  the image $\bar \C$. 
  The topological covering will be glued using  branched coverings  
  of a hemisphere by pants (we develop the 
  appropriate technique below). The proof of Theorem 2 relies on the 
  connectivity of the moduli spaces for the above branched coverings, 
  cp.~\cite {N2}. 
  
  \subheading 
  {3.1. On the space of branched covering of a hemisphere by a Riemann surface 
  with a boundary} 
  We start with some constructions. Denote by $\Lambda^+$ the upper 
hemisphere $\{z\in \bar \C \; \vert \;\text{Im}\; z \gs 0 \}$, and by
$P$ a  genus $g$ topological surface with a boundary consisting of $k$
connected components. Consider the set $ \H_{g,m}^k$ of all generic 
degree $m$ branched coverings of the form $f \: P\to \Lambda^+$. Let 
$a_{1},\ldots,a_k$ be all the distinct connected components of $\partial  
 P$. Given a partition $(m_1,...,m_k)\vdash m$,  
denote by $\H_{g,m}^k(m_{1},\ldots,m_k)\subset \H_{g,m}^k$ 
 the subset of maps 
$f\: P \to  \Lambda^+$ such that 
$\deg f\vert_{a_{i}}=m_{i}$ for $i=1,\ldots , k$. 
Obviously, 
$$
\H_{g,m}^k=\bigcup_{(m_1,...,m_k)\vdash m}\H_{g,m}^k(m_{1},
\ldots,m_k). 
$$

  Let $\widetilde P$ be a compact genus $g$ topological surface. Consider,  
in parallel, the set $\widetilde\H_{g,m}^k$ of all degree $m$  branched 
coverings $\tilde f\: \widetilde P\to \bar \C$ satisfying  the additional 
condition that  
all the ramification points are concentrated on 
$ \Lambda^+ \bigcup (-\frak {i})$, where $\frak i$ is the imaginary unit. 
Consider the obvious restriction map 
$\Psi\: \widetilde\H_{g,m}^k\to \H_{g,m}^k$, i.e. 
$\Psi(\tilde f\: \widetilde P\to \bar C)=(f\:P \to 
\Lambda^+)$, where $P= \tilde f^{-1}(\Lambda^+)$ and 
$f=\tilde f\vert_{P}$. According to \cite{N2, Theorem~4.1}, there exist  
unique complex structures on $P$ and $\tilde P$ for which the above
mentioned branching coverings are holomorphic. Thus we can consider
$\H_{g,m}^k$ and $\widetilde\H_{g,m}^k$ as spaces of meromorphic functions.

\proclaim {Lemma 1} 
The map $\Psi$ is a bijection{\rm ;} moreover, two maps $f_{1}$ and $f_{2}$ 
in $\H_{g,m}^k$ are equivalent if and only if their images 
$\Psi(f_{1})$ and $\Psi(f_{2})$ are equivalent. 
\endproclaim 

\demo {Proof} Denote $\Lambda^- =\{z\in \bar C\;\vert\; \text {Im}\; z 
\ls 0\}$ and fix a holomorphic degree $j$ map 
$\xi_{j}\: \Lambda^- \to \Lambda^-$ preserving $-\frak i$ and having no other 
ramification points on $\Lambda^-$ (such a $\xi_{j}$ obviously exists). 
 Take now an arbitrary function $f\in \H_{g,m}^k$.  It is 
always 
possible to identify each $a_{i}$ with $\partial \Lambda^-$ in such a 
way that $\xi_{m_{i}}\vert_{\partial \Lambda^-}=f\vert_{a_{i}}$. 
 Glueing copies of $\Lambda^-$ to all holes in $P$  gives a 
surface $\widetilde P$ without a 
 boundary.  At the same time, glueing  $f$ and $\xi_{m_{i}}$'s 
together  gives a new function $\tilde f\in \widetilde\H_{g,m}^k$. 
 Obviously, $\Psi(\tilde f)=f$, and moreover, this construction 
 sends equivalent functions to equivalent functions. 
\qed 
\enddemo

\proclaim {Lemma 2} 
 For any partition $(m_1,...,m_k) \vdash m$, the space  
$\H_{g,m}^k(m_{1},\ldots,m_k)$ is connected.
\endproclaim 

\demo {Proof} Define the corresponding set 
$\widetilde\H_{g,m}^k(m_{1},\ldots,m_k)\subset  \widetilde\H_{g,m}^k$ with 
$\tilde f^{-1}(-\frak i)$ consisting of a $k$-tuple of critical points 
$a_{1},\ldots,a_k$ with $\deg \tilde f\vert_{a_{i}}=m_{i}$. 
According to \cite {N1, N2}, the set $\widetilde\H_{g,m}^k(m_{1},\ldots,m_k)$ 
is connected. Therefore, by Lemma 1 one gets that the set  
$\H_{g,m}^k(m_{1},\ldots,m_k)\subset \H_{g,m}^k$ 
is connected as well. (The simplest case of this result is called 
the L\"uroth--Clebsch theorem, see e.g. \cite {Hu, Kl}.)
 \qed 
\enddemo

\proclaim {Lemma 3} With the above notation,
any $f \in \H_{g,m}^k(m_{1},\ldots,m_k)$ has exactly 
$m+k+2g-2$ simple ramification points on 
$\Lambda^+\setminus \partial\Lambda^+$. 
\endproclaim 
 
 \demo {Proof} Follows directly from the Riemann--Hurwitz formula 
 $$
\chi(P)+\sharp_{f}=
 \chi(\Lambda^+)\deg f,
$$
 where $\sharp_{f}$ is the number of simple ramification points of $f$, 
 $\deg f=m$ is the degree of $f$, $\chi(P)=2-2g-k$ 
 is the Euler characteristic of $P$, and $\chi(\Lambda^+)=1$ 
 is the Euler characteristic of $\Lambda^+.$  
 \qed 
\enddemo

\proclaim {Lemma 4} Given a genus $0$ surface $P$ with 
$k$ boundary components and a partition $(m_1,...,m_k)\vdash m\gs 3$,  
the Hurwitz number of $\H_{0,m}^k(m_{1},\ldots,m_k)$ 
equals 
$$
\frac {m^{k-3}(m+k-2)!m_1^{m_1}\dots m_k^{m_k}}
{m_1!\dots m_k!s(m_1,\dots,m_k)},
$$
where $s(m_1,\dots,m_k)$ is the number of symmetries of the set 
$\{m_1,\dots,m_k\}$.
\endproclaim 
 
 \demo {Proof} It follows immediately from Lemma 1 that the Hurwitz numbers
for the spaces $\widetilde\H_{0,m}^k(m_{1},\ldots,m_k)$ and 
 $\H_{0,m}^k(m_{1},\ldots,m_k)$ coincide. 
For the former space, this number equals  
 $\frac {m^{k-3}(m+k-2)!m_1^{m_1}\dots m_k^{m_k}}
{m_1!\dots m_k!s(m_1,\dots,m_k)}$ and was apparently known to Hurwitz,
 see \cite {GJ, St}.
 \qed 
\enddemo

{\smc Remark}. It is easy to see that for $m=2$ and $k=1,2$ the Hurwitz
number of $\H_{0,m}^k(m_1,\dots,m_k)$ equals 1, while the expression in
Lemma~4 gives 1/2. On the other hand, for $m=k=1$, the expression gives 1,
which is the correct answer.

\subheading
{3.2. Proof of Theorem~1}
Given a set $\Sigma$ of $2n-2$ distinct complex numbers of which exactly 
$2l$ are real, and a garden $G$ of order $2l$ with total weight $n$, 
we want to 
construct a topological branched covering $\bar \C\to\bar \C$ invariant 
under complex conjugation whose set of ramification points coincides with 
$\Sigma$, and whose garden is isomorphic to $G$. This will prove the Theorem,
since by \cite{N2, Theorem~4.1},
there exists a unique complex structure on $\bar \C$ for which this 
topological covering is holomorphic.

Consider $G$ as a represented garden, and
 let $\Pi(G)=\bigcup_{i=1}^q P_i$  
denote the weighted pants collection of $G$, see \S 2. In order to construct
a required topological branched covering $\bar \C\to\bar \C$,
we perform the following four steps.

{\smc Step 1.} Distribute $2l$ real numbers from $\Sigma_R$ between the
vertices of $G$, as described in the construction of $G(f)$ in \S 2.2.

{\smc Step 2.} Distribute $n-l-1$ complex conjugate pairs of numbers from 
$\Sigma_I$  between all pants in $\Pi(G)$. Let $P_i\in\Pi(G)$ be pants
with $c_i$ boundary components, and let $m_{ij}$ be the weight of the $j$th
boundary component. According to Lemma 3, we assign to $P_i$ exactly
$\sum_{j=1}^{c_i}m_{ij}+c_i-2$ pairs from $\Sigma_I$.

{\smc Step 3.} For any pants  $P_i\in\Pi(G)$ build a map
$f_i\:P_i\to\Lambda^{\pm}$ with prescribed ramification points. If $P_i$
are positive, then $f_i$ belongs to the space
$\H_{0, \mu_i}^{c_i}(m_{i,1},...,m_{i,c_i})$, where $\mu_i=\sum_j m_{ij}$; it
maps $P_i$ to $\Lambda^+$, and the ramification points 
are chosen as follows: from the each conjugate pair
assigned to $P_i$ on the previous step we take the point belonging to
$\Lambda^+$.  If $P_i$
are negative, then $f_i$ maps it to $\Lambda^-$, and the ramification
points are chosen in a similar way in $\Lambda^-$.

{\smc Step 4.} Glue all $f_i$'s together to get a map of the hemisphere 
containing $G$ to $\bar \C$ and, finally, glue the latter map  
   with its complex conjugate copy along the boundary of the hemisphere 
to get the actual branched covering  $\bar \C\to\bar \C$. 
    
 Let us explain the fourth step in detail.
Consider first  the case $l=0$. Taking the unique pants, called {\it basic},
 whose outer boundary is the circle of $G$ (identified with 
$\bar \R$ in the preimage), we glue to its map the maps of all 
 its neighboring pants  
by identifying these maps along their common boundary ovals.  Since by our 
construction the 
multiplicities of two maps having a common oval coincide on this oval,
the glueing process is possible (a similar procedure is 
used in the proof of Lemma 1). Having glued  the maps of all the neighbors 
to that of the basic pants, we continue with the neighbors of the neighbors, 
etc. 

 In the general case,  notice that each face $r$, $r=1,\dots,l+1$, 
contains the unique pants (called {\it basic} for $r$) 
 whose boundary coincides with that of the face $r$, see Fig~2. We can first 
  glue together the maps of all basic pants  and then continue as above. 
  The maps of a neighboring pair of basic 
  pants are glued together along their unique common arc which should be 
mapped to the prescribed  segment of $\bar \R$ in the image $\bar \C$ 
between the corresponding real ramification points, i.e. those 
  labeling the endpoints of the arc under consideration (the labels are
obtained on Step 1). 
Observe that these ramification points are regular points for each of
$f_i$'s, and become critical points only after glueing basic pants
together.
The weight of the arc defines the number of complete turns which this 
  arc should do around $\bar \R$ in the image, and its direction shows the 
orientation of the image of the arc. Thus the image of the arc is completely 
  determined by $G$. 
  
  Having glued all $f_i$'s together, we get a map $f$ from the 
  disc containing $G$ (identified with the upper hemisphere)
  to $\bar \C$. We take another copy of this disc (identified with the lower 
hemisphere) with the conjugate map 
  $\bar f$, and glue two hemispheres along $\bar \R$ into a sphere $\bar \C$ 
  with the final map $\bar \C\to \bar \C$ consisting of $f$  and $\bar f$. 
  
  One can easily see that the final map is the topological branched covering 
  with all properties  required by Theorem 1, and we are done. 

\subheading{3.3. Connected components of $\R \H_{0,n}$}
 The following statement is crucial for the proof of Theorem 2. 
 
 \proclaim {Lemma 5} Two generic real rational functions  
 $f_i\: \bar \C\to \bar \C$, $i=0,1$,  are equivalent 
  if and only if 
 
 a) their gardens $G(f_0)$ and $G(f_1)$ are isomorphic, i.e. there exists
a bijection
of the vertex sets of $G(f_0)$ and $G(f_1)$ that preserves
chords, their orientation, labels, rooted trees,
and their weights;

   b) the restrictions of $f_0$ and $f_1$ to each pair of pants
  identified by the above isomorphism of gardens are equivalent. In particular,
   the sets of complex critical values assigned to each pair of pants
  identified by the above isomorphism of gardens coincide.
   \endproclaim

 \demo {Proof} Obviously, if $f_0$ and $f_1$ are equivalent 
 then a)-b) are automatically satisfied. On the other hand, using condition  
 b) we can construct, for each pair of pants identified by the above 
 isomorphism
of gardens, a homeomorphism making the restrictions of $f_0$ and $f_1$ 
to these pants equivalent.  Now using a) we can glue together 
these homeomorphisms defined 
on pairs of pants into a global homeomorphism $\bar \C\to \bar \C$ making 
$f_0$ and $f_1$ equivalent. As usual, the constructed
homeomorphism provides a biholomorphic map by inducing the complex
structure on the preimage $\bar \C$ from that on the image $\bar
\C$.
\qed 

\enddemo 
\medskip
\demo{Proof of Theorem 2} 
Let us show the easy implication first. Assume that two 
rational functions $f_0$ and $f_1$ belong to the same connected component 
of $\R \H_{0,n}$. Let us show that $G(f_0)$ is equivalent to $G(f_1)$. 
Take some path $f_t\subset \R \H_{0,n}$, $t\in[0,1]$, 
connecting $f_0$ and $f_1$. The only data related to the gardens 
of real rational functions which can vary along 
$f_t$ (up to a diffeomorphism of $\bar \C$ invariant under complex 
conjugation) are the values of ramification points. But since they never 
collide, one gets that the complex ramification points remain complex, the   
 real ramification points remain real, and can only experience a cyclic shift. 
 Thus, $G(f_0)$ is equivalent to $G(f_1)$.  
   
Conversely, take two functions $f_0$ and $f_1$ in $\R \H_{0,n}$ 
whose gardens $G_0=G(f_0)$ and $G_1=G(f_1)$ are equivalent. 
Notice that the equivalence of $G_0$ 
and $G_1$ implies the 1-1-correspondence between the sets $\Pi(G_0)$ and 
$\Pi(G_1)$ of the weighted pants collections, i.e. the existence of  a 
1-1-correspondence between pants for $f_0$ and $f_1$. 

Let $\Sigma_0$ and $\Sigma_1$ denote the sets of ramification points 
of $f_0$ and $f_1$, respectively.

The equivalence of $G_0$ and $G_1$ implies that $\Sigma_0$ and $\Sigma_1$
belong to the same connected 
component of $\R Q_{0,n}$. Moreover, we can connect $\Sigma_0$ and $\Sigma_1$ 
by a path 
$\Sigma_t$, $t=[0,1]$, in this component in such a way that 
for any $i$, the subset of $\Sigma_0$ corresponding to the $i$th pants
in $\Pi(G_0)$ will be transformed along $\Sigma_t$ into the 
subset of $\Sigma_1$ corresponding to the $i$th pants
in $\Pi(G_1)$.
 Using the covering homotopy property of $\R\Phi_{0,n}$ (see \cite{N2}) 
over the path $\Sigma_t$, 
we get another rational map $\tilde f_1$ which lies in the same 
connected component of $\R \H_{0,n}$ as $f_0$; therefore, the garden 
$\widetilde G_1$ 
of $\tilde f_1$ is equivalent to $G_0$ (and hence to $G_1$)
by the first part of this proof. Moreover, the set of ramification points
of $\tilde f_1$ coincides with $\Sigma_1$, and the distributions of
ramification points among pants for $f_1$ and $\tilde f_1$ are identical. 

Now for each pants from $P_i(\widetilde G_1)$ we can, using the connectivity 
of the space of maps proved in Lemma 2, find a path between the restriction 
of $\tilde f_1$ to these pants
and the restriction of $f_1$ to the 
corresponding pants from $\Pi(G_1)$ that keeps the restrictions of 
$f_1$ and $\tilde f_1$  to all other pants unchanged. Doing this procedure for
  every pants we connect $\tilde f_1$ with a map $\bar f_1$, which  
together with $f_1$ satisfies all the conditions of Lemma 5. Therefore,
$\bar f_1$ is equivalent to $f_1$, and is connected with $f_0$ by a path in
$\R\H_{0,n}$, hence $f_0$ and $f_1$ belong to the same connected component
of $\R\H_{0,n}$.
\qed 
\enddemo 

\vglue 12pt 
\line{\bf 4. Hurwitz Numbers \hfil}
\vglue 5pt

To find out the number of nonequivalent functions corresponding to the 
same garden, consider a weighted rooted tree $T$. The node set of $T$
is $\{0,1,\dots,k\}$ for some $k\gs0$, and $0$ is the root of $T$. The weight 
of node $i$ is denoted $w_i$. Let $i$ be an arbitrary node, and $i_1,
\dots,i_c$ be all of its children. The number of children is called
the {\it degree\/} of the node $i$ and is denoted $c_i$ (see Figure~3).

\vskip 10pt
\centerline{\hbox{\epsfxsize=5cm\epsfbox{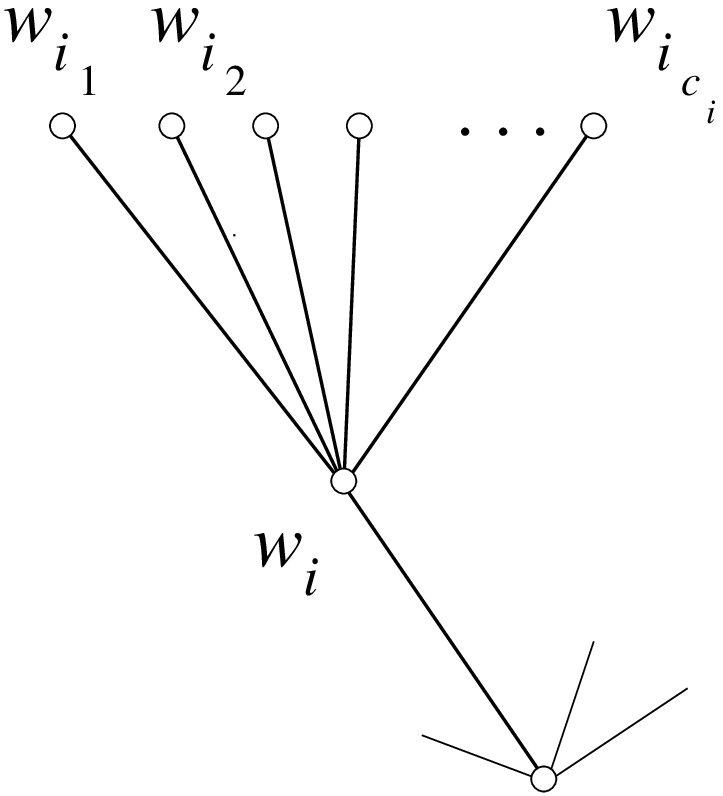}}}
\midspace{1mm}
\caption{Fig.~3. A node and its children}
\vskip 5pt 

Besides, we define the
{\it modified weight\/} $\widetilde w_i$ as the sum of the weight of $i$ 
and the weights of all of its children. The {\it total weight\/}
$w_T$ of the tree $T$ is equal to the sum of the modified weights of all
of the nodes, including the root. 
Finally, we define the {\it symmetry
factor\/} $s_i$ as the number of automorphisms of
the set $\{w_i, w_{i_1},\dots,
w_{i_{c_i}}\}$.  The {\it Hurwitz number\/} $H_T$ of the tree $T$ is defined
by
$$
H_T=\frac{(w_T-1)!w_0!2^{e_T}}{w_0^{w_0}}
\prod_{i=0}^k\frac{w_i^{2w_i}\widetilde w_i^{c_i-2}}{(w_i!)^2s_i},
$$
where $e_T$ is the number of nodes in $T$ whose modified weight equals 
$2$.

Assume first that the set $\Sigma$ does not contain real numbers, thus 
$l(\Sigma)=0$.
The garden $G(f)$ for an arbitrary function $f$ such that 
$\Sigma(f)=\Sigma$ has order $0$
and consists of a trivial chordless chord diagram and a single tree.
Such a garden we call an {\it imaginary\/} garden. Observe that the total
weight of an imaginary garden equals the total weight of its single
tree.

\proclaim{Theorem 6}  Let $\Sigma$ be an arbitrary set of $2n-2$ distinct
complex numbers invariant under the standard involution and containing no 
real numbers.  
Let $G$ be an imaginary garden of total weight $n$
and $T$ be its
single tree. The number of topologically nonequivalent functions
$f\in\R\H_{0,n}$ such that $\Sigma(f)=\Sigma$ and $G(f)=G$ is equal to the 
Hurwitz number $H_T$.
\endproclaim

\demo{Proof} By Lemma 5, we have to calculate the number of ways to
execute Steps 2 and 3 in the proof of Theorem 1. First, we distribute 
the elements of $\Sigma$ over the pants defined by $T$. By Lemma 
3, the number of ramification points corresponding to the pants $P_i$
equals $\widetilde w_i+c_i-1$. The total number of points to be
distributed equals
$$
\sum_{i=0}^k(\widetilde w_i+c_i-1)=w_{T}+\sum_{i=0}^k c_i -k-1=n-1,
$$
since the total weight of $T$ equals $n$ and $\sum_{i=0}^k c_i$ is
the number of inner nodes in $T$. Therefore,
the total number of the ways to execute Step 2 equals
$$
\frac{(n-1)!}{\prod_{i=0}^k (\widetilde w_i+c_i-1)}.
$$

The number of ways to execute Step 3 is described in Lemma 4 and the 
Remark following the lemma. Therefore, the total number in question
equals
$$
\frac{(n-1)!}{\prod_{i=0}^k (\widetilde w_i+c_i-1)}
\prod_{i=0}^k
\frac{\widetilde w_i^{c_i-2}(\widetilde 
w_i+c_i-1)!w_i^{w_i}w_{i_1}^{w_{i_1}}\dots w_{i_c}^{w_{i_c}}}
{s_i w_i!w_{i_1}!\dots w_{i_{c_i}}!}\cdot 2^{e_{T}}.
$$
The expression $w_i^{w_i}/w_i!$ for a given inner node $i$ enters the 
latter product twice: once when the node itself is considered, and once 
more when the node appears as a child of its parent node. After 
cancellations we get the desired result.
\qed
\enddemo

In the general case, let $G$ be a garden of order $2l$ and total weight 
$w$, 
and $T_1,\dots,T_{l+1}$ be its trees. The Hurwitz number of the 
garden is defined by
$$\align
H_G&=(w-l-1)!\prod_{i=1}^{l+1}\frac{H_{T_i}}{(w_{T_i}-1)!}\\&=
2^{e_G}(w-l-1)!\prod_{i\in R_G}\frac{w_i!}{w_i^{w_i}}
\prod_{i\in N_G}\frac{w_i^{2w_i}\widetilde w_i^{c_i-2}}{(w_i!)^2s_i},
\endalign
$$
where  $N_G$ is the set of the nodes of the trees $T_1,\dots,T_{l+1}$, 
$R_G$ is the set of the roots of these trees, and
$e_G$ is the number of nodes in $N_G$ whose modified weight equals $2$.

The following proposition follows easily from Lemma 5 similarly to
Theorem 6.

\proclaim{Theorem 7} Let $\Sigma$ be an arbitrary set of $2n-2$ distinct
complex numbers invariant under the standard involution, of which 
exactly $2l(\Sigma)$ are real, and let $G$ be a garden of order $2l(\Sigma)$ 
and total weight $n$. The number of topologically nonequivalent functions
$f\in\R \H_{0,n}$ such that $\Sigma(f)=\Sigma$ and $G(f)=G$ is equal to the 
Hurwitz number $H_G$. 
\endproclaim
\vglue 5pt
\line{\bf 5. Final Remarks \hfil}
\vglue 5pt

In this note we assigned to each connected component in the space 
$R \H_{0,n}$ a combinatorial object called a garden. Unfortunately, to count 
the total number of all gardens of a given weight $n$ seems to be a difficult 
problem. Two cases look more accessible, namely, the {\it elliptic} case when 
no critical values are real, and the {\it hyperbolic} case when all critical 
values  are real.

{\smc Problem.} Count the number of connected components in the space of all 
hyperbolic and elliptic generic functions of degree $n$. 

In the hyperbolic case, the major combinatorial difficulty is to count the 
total number of admissible labelings of a given planar chord diagram with  
$2n-2$ vertices, see \cite {SV2}. In the elliptic case, one should count the
total number of nonisomorphic planar trees with total weight $n$.

\end